\documentclass[12pt]{article}
\usepackage[colorlinks]{hyperref}
\usepackage{color}
\usepackage{graphicx}
\usepackage{graphics}
\usepackage{makeidx}
\usepackage{showidx}
\usepackage{latexsym}
\usepackage{amssymb}
\usepackage{verbatim}
\usepackage{amsmath}
\usepackage{amsthm}
\usepackage{amsfonts}
\usepackage{geometry}

\title{Collision of a Hard Ball with Singular Points of the Boundary}
\author{H. Attarchi and L.A. Bunimovich}
\begin{document}
	\maketitle
	\noindent
\begin{abstract}
	Recently were introduced physical billiards where a moving particle is a hard sphere rather than a point as in standard mathematical billiards. It has been shown that in the same billiard tables the physical billiards may have totally different dynamics than mathematical billiards. This difference appears if the boundary of a billiard table has visible singularities (internal corners if the billiard table is two-dimensional), i.e. the particle may collide with these singular points. Here, we consider the collision of a hard ball with a visible singular point and demonstrate that the motion of the smooth ball after collision with a visible singular point is indeed the one that was used in the studies of physical billiards. So such collision is equivalent to the elastic reflection of hard ball's center off a sphere with the center at the singular point and the same radius as the radius of the moving particle.\par
	\noindent{\bf Keywords:}
	physical billiards, mathematical billiards, visible and invisible singularities
\end{abstract}
%%%%%%%%%%%%%%%%%%%%%%%%%%%%%%%%%%%%%%%%%%%%%%%%%%%%%%%%%%%%%%%%%%%%%%%%%%%%%%%%%%%%%%%%%%%%%%%%%%%%%%%%%%%%%%%
\section{Introduction}
Mathematical billiards serve as relevant models of various phenomena in mechanics, geometric optics and acoustics, statistical physics, and quantum physics. Such billiards also constitute one of the most popular and arguably the most visual class of dynamical systems in the mathematical studies. In mathematical billiards, a point particle moves by inertia in a domain with boundary. When the point particle reaches the boundary, it gets elastically reflected.\par
Recently were introduced physical billiards where the moving particle is a hard ball \cite{Bun19}. It was shown in this paper that in the transition from a mathematical to a physical billiard in the same billiard table any type of transition from chaotic to regular dynamics and vice versa may occur. Moreover, such transition from the point to a finite size particle can completely change the dynamics of some classical and well-studied models like e.g. the Ehrenfests' Wind-Tree model \cite{ABB}. In quantum systems, a ``particle'' naturally has a finite size due to the uncertainty principle which leads to some new findings in the quantum chaos theory \cite{PCB,RBG}.\par
Interesting changes in dynamics occur when the boundary of a billiard table has a visible singularity, i.e. a point in the intersection of two or more smooth components of the boundary such that a small enough physical particle can hit that point of the boundary. If a billiard table is two-dimensional, then such singularities are internal corners where two smooth components of the boundary intersect and make an angle greater than $\pi$ inside the billiard table. In all papers cited above, it was assumed that reflection of the ball off such visible singularity occurs in a natural manner corresponding to the simplest elastic collision. In the present note, we justify this assumption for a smooth hard ball. It is worthwhile to mention that there are other types of reflection of a ball off a visible singular point that correspond to a rough ball which may acquire rotation after such collision \cite{Gar69} even under the assumption that it is a no-slip collision \cite{BG93,CF17}.
%%%%%%%%%%%%%%%%%%%%%%%%%%%%%%%%%%%%%%%%%%%%%%%%%%%%%%%%%%%%%%%%%%%%%%%%%%%%%%%%%%%%%%%%%%%%%%%%%%%%%%%%%%%%%%%
\section{Different types of boundary singularities in billiard tables}
Let $Q$ be a domain in $d$-dimensional Euclidean space $\mathbf{R}^d$ such that its boundary $\partial Q$ is the union of a finite number of $C^1$-smooth $(d-1)$-dimensional manifolds. A point $q$ of the boundary $\partial Q$ is called singular if the boundary is not differentiable at that point. That means a singular point belongs to the intersection of some (at least two) differentiable (aka regular) components of the boundary. Note that we also call a singular point in dimension two (i.e. $\dim Q=2$) a corner. All non-singular points of the boundary $\partial Q$ are called regular points.\par
Consider a free motion of a hard ball (a disk in $\dim 2$) of radius $r>0$ in the domain $Q$ with elastic reflections off the boundary $\partial Q$. The resulting dynamical system is called a physical billiard \cite{Bun19}, and the domain $Q$ is a billiard table. To describe the dynamics of such ball, it is enough to follow the motion of its center. It is easy to see that the center of the ball moves in the smaller billiard table, which one gets by moving any point $q$ of the boundary by $r$ to the interior of the billiard table along the internal unit normal vector $n(q)$ \cite{Bun19}.\par
We will call a singular point $q$ of the boundary $\partial Q$ an invisible singular point if for any $r>0$ the hard ball of radius $r$ cannot hit that point. Otherwise, a singular point is called a visible singular point. Therefore, $q$ is a visible singular point of a billiard table if a ball with a sufficiently small radius can hit $q$. A formal mathematical definition of a visible singular point (in any dimension) is the following one. A singular point $A$ is a visible singular point if for any neighborhood $N$ of $A$ the convex hull of $Q\cap N$ contains a neighborhood of $A$.\par
\begin{figure}[ht]
	\begin{center}
		\includegraphics[width=7cm]{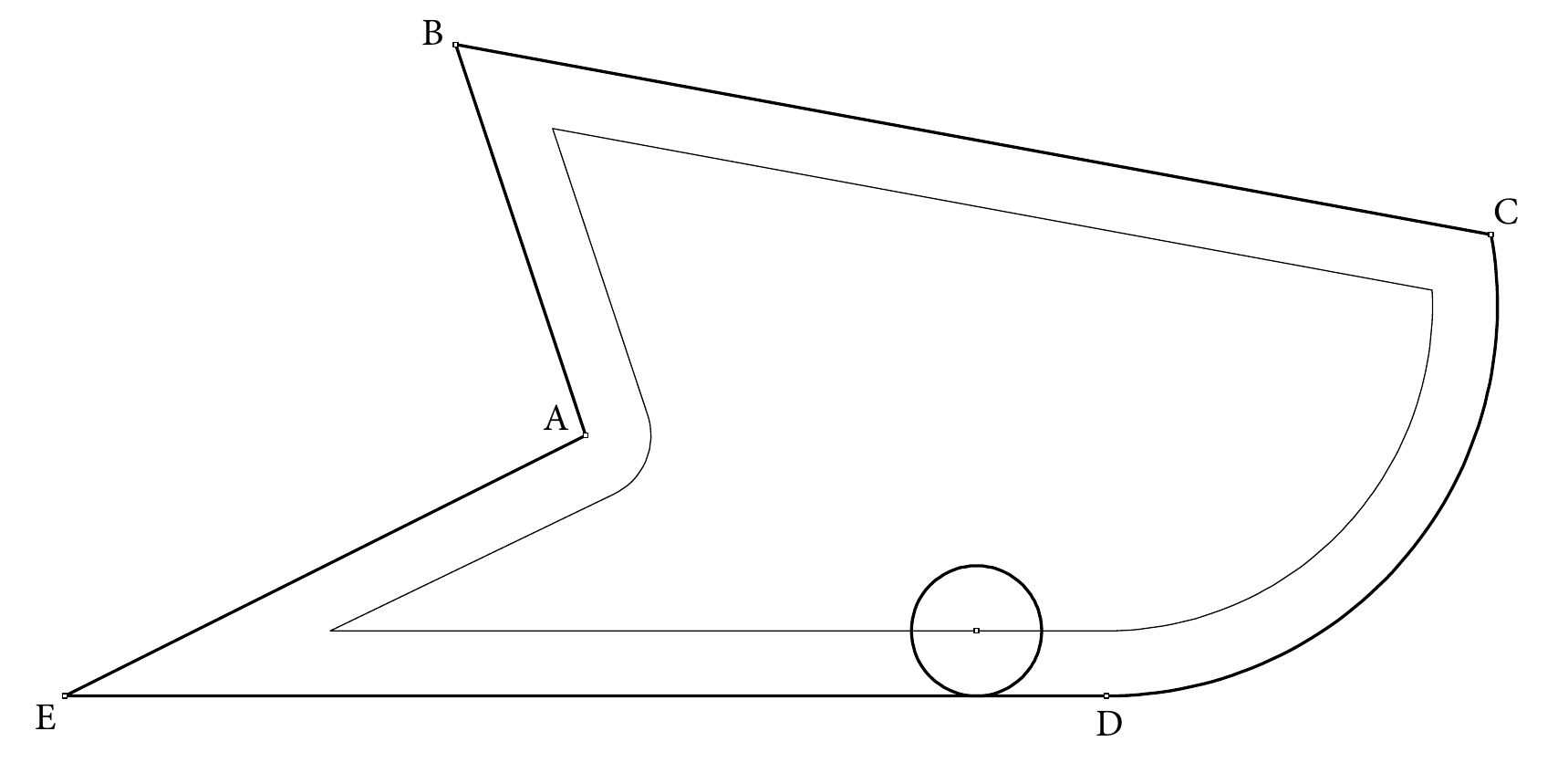}
		\caption{Corners (singular points) B, C, and E are invisible to any disk. The point D is not singular, since boundary is differentiable at D. The corner A is an internal corner (a visible singular point).}~\label{Corners}
	\end{center}
\end{figure}

We also call a visible singular point in dimension two an internal corner. For example, Fig. \ref{Corners} shows visible and invisible singular points in dimension two.\par
\begin{figure}[ht]
	\begin{center}
		\includegraphics[width=4cm]{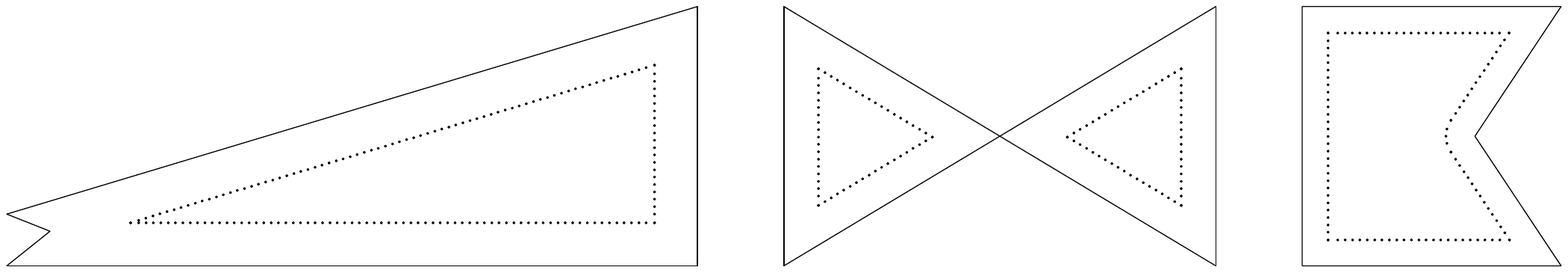}
		\caption{An internal corner becomes invisible when radius of disk is larger than some constant.}~\label{invis}
	\end{center}
\end{figure}

Note that being a visible singular point (an internal corner) does not mean that a hard ball of any radius $r>0$ can reach (hit) that point. Namely, if the radius of the particle is larger than some constant (which depends on the shape of a billiard table), then some visible singular points become invisible (see Fig. \ref{invis}).\par
\begin{figure}[ht]
	\begin{center}
		\includegraphics[width=7cm]{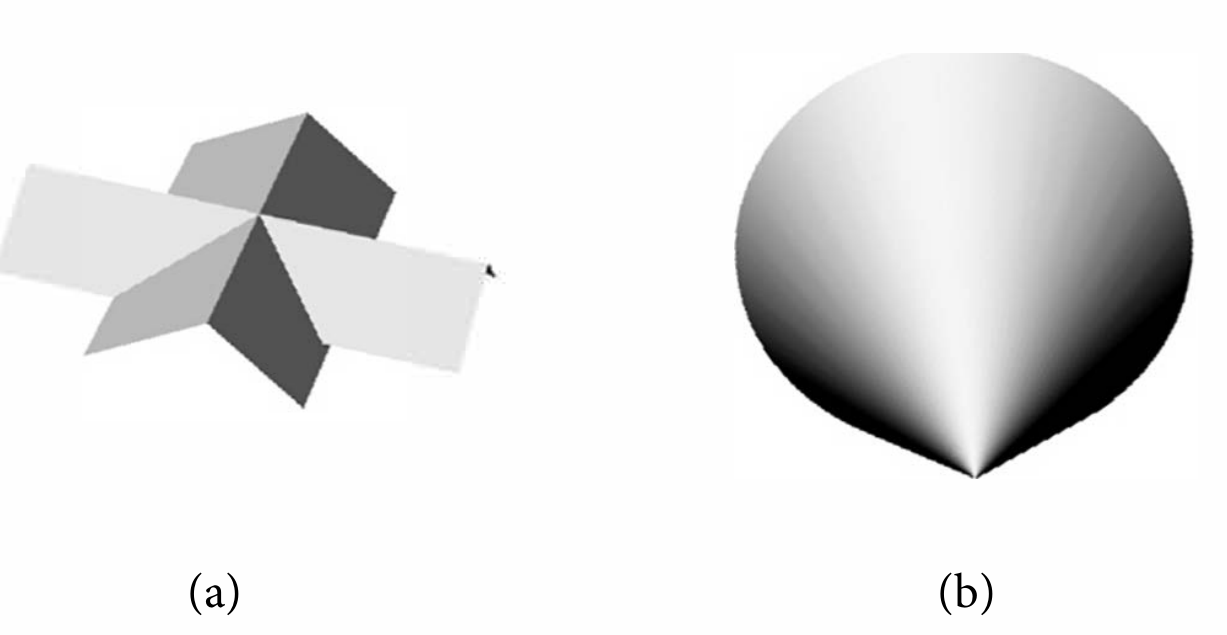}
		\caption{(a) There are two lines of visible singular points. When a hard ball hits a point on those lines, its center is on an arc of a circle centered at that point and orthogonal to the corresponding line. However, at the moment of collision with the intersection point of those two lines the center of  hard ball can be only in one position. (b) Here is one isolated visible singular point. At the moment of collision with such singularity the center of a hard ball is on a piece of $2$-sphere centered at that singular point.}
		\label{3dim}
	\end{center}
\end{figure}

Observe that at the moment of collision with a visible singular point, the center of hard ball can be at different positions, and these possible positions depend on the shape of the boundary $\partial Q$ (see Fig. \ref{3dim}). This should be contrasted with the collision of the ball off the boundary at a regular point, when the center of the ball always has one position, namely at the distance $r$ on the internal normal line to the boundary of a billiard table. In Fig. \ref{3dim}, two situations are depicted, which may happen in three dimensional billiard tables.\par
Since the particle is a hard ball, it will keep its shape at the moment of collision. Hence the center of the hard ball is at the distance $r$ from a collision point (regardless of whether this point is a regular or singular point of the boundary). Therefore, the boundary of the reduced billiard table of the mathematical billiard, which has the same dynamics as the considered physical billiard \cite{Bun19}, acquires a piece of a sphere (or an arc of a circle if the dimension of the billiard table is two) of radius $r$ with the center at the visible singular point. Hence the reduced billiard table of the equivalent mathematical billiard has a dispersing component in the boundary which generates a chaotic (hyperbolic) dynamics in case if a moving particle is a smooth hard ball.\par
However, for any type of a hard ball, a reduced billiard table of the equivalent mathematical billiard acquires a dispersing (or semi-dispersing) component. This fact holds true for any type of collision of the physical ($r>0$) particle with the boundary at a visible singular point. However, such collisions can generally be elastic or inelastic and with or without slip \cite{CF17,Gar69}. Dynamics of rough ball even in case of no-slip collisions is much more complicated than the dynamics of a smooth ball.\par
\begin{figure}[ht]
	\begin{center}
		\includegraphics[width=7cm]{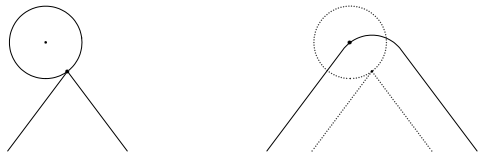}
		\caption{A collision between a disk and a visible singular point (here, an internal corner) is shown in the left picture. On the right, one can see its equivalent for a (virtual) collision between disk's center and an arc of a circle centered at the internal corner with the same radius as the disk's radius.}~\label{internalcorner}
	\end{center}
\end{figure}

In Fig. \ref{internalcorner}, it is easy to see the boundary of the reduced billiard table of the mathematical billiard acquires a dispersing component, because of the case of dimension two depicted. Here, the center of a disk can be located at any point of an arc of the circle with the center at the singular point and with the radius equals to the radius of the disk.
%%%%%%%%%%%%%%%%%%%%%%%%%%%%%%%%%%%%%%%%%%%%%%%%%%%%%%%%%%%%%%%%%%%%%%%%%%%%%%%%%%%%%%%%%%%%%%%%%%%%%%%%%%%%%%%
\section{No-slip collisions of a hard ball with a visible singular point}
In the case of no-slip collisions, each reflection of the moving particle (hard ball) off the boundary occurs at a single point. Hence a collision at any point of the boundary does not depend on the shape of the boundary elsewhere. Therefore, the collision problem can be actually considered as a reflection of a hard ball off a point.\par
At the moments of the collision, the impulse $\Delta P$ decomposes into two components, which are the normal impulse $\Delta P_N$ acting towards the center of hard ball and the tangent impulse $\Delta P_T$ based on friction which is tangent to the hard ball at the collision point. The tangent impulse can result in either loss of kinetic energy or exchange between linear and angular momentum while the total kinetic energy is preserved. We will consider the friction-free (elastic) collision and the case when the impulse $\Delta P_T$ results in an exchange between linear and angular momentum without loss of energy. In other words, we consider only conservative (Hamiltonian) dynamics.\par
Let a hard ball of radius $r>0$ with the center at a point $O$ and mass $m=1$ hits a visible singular point $A$ of the boundary of a billiard table $Q$. Denote the linear velocity of hard ball's center just before (after) the collision by $V^b$ ($V^a$). Consider now a decomposition of $V^b$ to two components $V_N^b$ and $V_T^b$, where $V_N^b=Proj_{\overrightarrow{OA}}V^b$ and $V_T^b=V^b-V_N^b$. Note that we will use the superscript $a$ instead of $b$ to denote velocity components at a moment of time right after the reflection. Denote also the vector form of angular velocity just before (after) the collision about the point $O$ by $\omega^b$ ($\omega^a$).\par
The collision map $S$ at point $A$ will map linear components and the angular component of the velocity just before collision $(V_N^b,V_T^b,\omega^b)$ to those right after collision $(V_N^a,V_T^a,\omega^a)$. The map $S$ has the following properties:
\begin{enumerate}
	\item The map $S$ is an orthogonal map because of the assumption that the system in question is Hamiltonian.
	\item Because of time reversibility of dynamics, $S^2$ is the identity map.
	\item The normal component of the linear velocity with respect to the boundary of hard ball at the contact point $A$ (i.e. $V_N^b$) always reverses under the map $S$.
\end{enumerate}
The conditions 1 and 2 imply that the eigenvalues of the map $S$ are $1$ or $-1$. In view of 3, one gets $S(V_N^b,V_T^b,\omega^b)=(-V_N^b,V_T^a,\omega^a)$, or equivalently, $V=(V_N^b,\vec{0},\vec{0})$ is an eigenvector of $S$ corresponding to the eigenvalue $-1$. It also implies that $\Delta P_N=-2V_N^b$. Without any loss of generality we assumed that the mass of a hard ball is $1$.\par
The Hamiltonian system under consideration satisfies three conservation laws of the kinetic energy $K$, the linear momentum $P$, and the angular momentum $L$ about the point $O$. These conservation laws in dimension $3$ are given by the relations
\begin{equation}\label{eq main}
	\left\{
	\begin{array}{l}
		K^b=\frac{1}{2}\left(|V_N^b|^2+|V_T^b|^2+I|\omega^b|^2\right)
		\\ \hspace{0.5cm}=\frac{1}{2}\left(|V_N^a|^2+|V_T^a|^2+I|\omega^a|^2\right)=K^a,\\
		P^b+\Delta P=V_N^b+V_T^b+\Delta P_N+\Delta P_T=V_N^a+V_T^a=P^a
		\\
		L^b+\Delta P_T\times\overrightarrow{AO}=I\omega^b+\Delta P_T\times\overrightarrow{AO}=I\omega^a=L^a,
	\end{array}\right.
\end{equation}
where $I$ is the moment of inertia of the hard ball.\par
Using that $V_N^a=-V_N^b$ and $\Delta P_N=-2V_N^b$, one can simplify (\ref{eq main}) as
\begin{equation}\label{main}
	\left\{
	\begin{array}{l}
		|V_T^b|^2+I|\omega^b|^2=|V_T^a|^2+I|\omega^a|^2,\\
		V_T^b+\Delta P_T=V_T^a
		\\
		I\omega^b+\Delta P_T\times\overrightarrow{AO}=I\omega^a.
	\end{array}\right.
\end{equation}
By solving (\ref{main}) for $\Delta P_T$, we get 
\begin{equation}\label{delta}
	\langle\Delta P_T,\frac{r^2+I}{I}\Delta P_T+ 2V_T^b+2\overrightarrow{AO}\times\omega^b\rangle=0,    
\end{equation}
where $\langle.,.\rangle$ is the inner product in $\mathbf{R}^3$ and $r$ is radius of hard ball. It is easy to see that $\Delta P_T=\vec{0}$ is a solution of (\ref{delta}) under the condition that there is no friction.\par
Observe that the conservation laws in dimension $2$ are the same as in (\ref{eq main}) under the assumption that the billiard table $Q$ is a subset of $xy$-plane in $\mathbf{R}^3$.
%%%%%%%%%%%%%%%%%%%%%%%%%%%%%%%%%%%%
\subsection{Friction-free collision (a smooth ball)}
In this section, we study a friction-free (i.e. $\Delta P_T=\vec{0}$) Hamiltonian system. In this case, (\ref{main}) implies
$$V_T^a=V_T^b,\hspace{1cm} \omega^a=\omega^b.$$
Here, the solution $(V_N^a,V_T^a,\omega^a)=(-V_N^b,V_T^b,\omega^b)$ of (\ref{eq main}) corresponds to the case of smooth hard ball \cite{Gar69} when the ball does not acquire rotation upon collision. Thus, in this case, we have an elastic reflection where the angle of incidence is equal to the angle of reflection.\par
Also, this friction-free collision is equivalent to the elastic reflection of the hard ball's center $O$ off a piece of a $2$-sphere (it can be an arc of a circle) centered at the visible singular point $A$ with the same radius as the radius of the hard ball \cite{ABB,Bun19}.\par
In case of dimension $3$, the collision map $S$ is a linear map from a $6$-dimensional subspace of $\mathbf{R}^9$ to itself with eigenvalues $1$ and $-1$. When $\Delta P_T=\vec{0}$, the eigenvectors which correspond to these eigenvalues have the forms $(\vec{0},V_T^b,\omega^b)$ and $(c\overrightarrow{AO},\vec{0},\vec{0})$, respectively, where $c$ is a constant. Also, the eigenspaces corresponding to the eigenvalues $1$ and $-1$ have dimensions five and one, respectively.
%%%%%%%%%%%%%%%%%%%%%%%%%%%%%%%%%%%%
\subsection{Collisions with friction (a rough ball)}
For the Hamiltonian system under consideration, the presence of the frictional force means that $|\Delta P_T|\neq0$. The corresponding solution of (\ref{eq main}) when $|\Delta P_T|\neq0$ describes dynamics of a rough ball \cite{Gar69}, which has ultra-elastic (no-slip) reflections off the boundary. In this case, the tangential component of the linear velocity partially transfers to the angular velocity and vice versa.\par
A nontrivial solution for $\Delta P_T$ in (\ref{delta}), is given by
\begin{equation}\label{rough}
	\Delta P_T=-\frac{2I}{r^2+I}(V_T^b+\overrightarrow{AO}\times\omega^b).
\end{equation}
Let $S$ be the collision map in dimension $3$ when the tangent impulse $\Delta P_T$ is given by (\ref{rough}). Then $(\vec{0},V_T^b,\omega^b)$ is an eigenvector of the collision map $S$ corresponding to the eigenvalue $1$ if $V_T^b+\overrightarrow{AO}\times\omega^b=\vec{0}$. The solution set of the vector equation $V_T^b+\overrightarrow{AO}\times\omega^b=\vec{0}$ is a three dimensional space. Hence, the eigenspace corresponding to the eigenvalue $1$ of the collision map $S$ is a $3$-dimensional space.\par
Moreover, $(\vec{0},V_T^b,\omega^b)$ is an eigenvector of the collision map $S$ which corresponds to the eigenvalue $-1$ if $V_T^b\times\overrightarrow{AO}-I\omega^b=\vec{0}$. In this case, the solution set of the vector equation $V_T^b\times\overrightarrow{AO}-I\omega^b=\vec{0}$ is a two dimensional space. This implies that the eigenspace corresponding to the eigenvalue $-1$ of the collision map $S$ is a $3$-dimensional space (we know $(\overrightarrow{AO},\vec{0},\vec{0})$ is another eigenvector for eigenvalue $-1$).
%%%%%%%%%%%%%%%%%%%%%%%%%%%%%%%%%%%%%%%%%%%%%%%%%%%%%%%%%%%%%%%%%%%%%%%%%%%%%%%%%%%%%%%%%%%%%%%%%%%%%%%%%%%%%%%
\section*{Acknowledgements}
The authors are grateful to R. Feres for helpful discussions.
%%%%%%%%%%%%%%%%%%%%%%%%%%%%%%%%%%%%%%%%%%%%%%%%%%%%%%%%%%%%%%%%%%%%%%%%%%%%%%%%%%%%%%%%%%%%%%%%%%%%%%%%%%%%%%%

\end{document}